\documentclass[12pt,psfig]{article}
\usepackage{graphicx,epsfig}
\usepackage{amsmath}
\usepackage{mathrsfs}
\usepackage{amssymb}
\def\ge{\geq}

\def\bbb{\begin{eqnarray*}}

\def\eee{\end{eqnarray*}}

\begin{document}

\baselineskip=17pt
\begin{center}

\vspace{-0.6in} {\large \bf On distributional chaos in
non-autonomous discrete systems}\\ [0.2in]

Hua Shao, Yuming Shi$^{*}$, Hao Zhu

\vspace{0.15in} Department of Mathematics, Shandong University \\
 Jinan, Shandong 250100, P.~R. China\\

\footnote{$^{*}$ The corresponding author.}
\footnote{$^{**}$ Email addresses: huashaosdu@163.com (H. Shao), ymshi@sdu.edu.cn (Y. Shi), haozhu@mail.sdu.edu.\\cn  (H. Zhu).} \

\end{center}

{\bf Abstract.} This paper studies distributional chaos in non-autonomous discrete systems
generated by given sequences of maps in metric spaces. In the case that the metric space is
compact, it is shown that a system is Li-Yorke $\delta$-chaotic if and only if it is distributionally
$\delta'$-chaotic in a sequence; and three criteria of distributional $\delta$-chaos are established,
which are caused by topologically weak mixing, asymptotic average shadowing property, and some expanding
condition, respectively, where $\delta$ and $\delta'$ are positive constants. In a general case, a criterion
of distributional chaos in a sequence induced by a Xiong chaotic set is established.\medskip

{\bf \it Keywords}:\ non-autonomous discrete system; distributional chaos; Li-Yorke chaos;
topologically weak mixing; shadowing property.
\medskip

{2010 {\bf \it Mathematics Subject Classification}}: 37B55, 37D45, 37A25.

\bigskip

\noindent{\bf 1. Introduction}\bigskip

In the present paper, we study the following non-autonomous discrete system (simply, NDS):
\vspace{-0.2cm}$$x_{n+1} = f_n(x_n), \;\;n \geq0,  \ \eqno(1.1)\vspace{-0.2cm}$$
where $f_n: X\to X$ is a map, and $(X,d)$ is a metric space with metric $d$.

When $f_n=f$ for each $n\ge 0$, (1.1) is the following autonomous discrete system (simply, ADS):
\vspace{-0.2cm}$$x_{n+1}=f(x_n), \;\;n \geq0,  \ \eqno(1.2)\vspace{-0.2cm}$$
where $f: X\to X$ is a map.

Note that the ADS (1.2) is governed by the single map $f$ while the NDS (1.1) is generated by iteration of a sequence of maps
$\{f_n\}_{n=0}^{\infty}$ in a certain order. Thus, it is more difficult to study dynamical behaviors of NDSs than those
of ADSs in general. However, many complex systems occurring in the real world problems such as physical, biological,
and economical problems are necessarily described by NDSs. For example, the well-known logistic system
\vspace{-0.2cm}$$x_{n+1}=rx_n(1-x_n), \;\;n \geq0,\vspace{-0.2cm}$$
describes the population growth under certain conditions. Note that
the parameter $r$ will vary with time when the natural environment changes.
So it is more reasonable to describe the population growth in this case by the following
non-autonomous system:
\vspace{-0.2cm}$$x_{n+1}=r_nx_n(1-x_n), \;\;n \geq0.\vspace{-0.2cm}$$
Hence, it is meaningful to study NDSs and many mathematicians focused on complexity of NDSs
in recent years [4, 5, 8, 12, 18, 19, 24, 25, 27, 31--33, 38, 44].

Li and Yorke firstly introduced the concept of chaos in their famous work ``period three implies chaos" [13],
which has activated sustained interest in the frontier research on discrete chaos theory. Later, various definitions
of chaos were developed, such as Devaney chaos [7], Auslander-Yorke chaos [2], generic chaos [29], dense chaos [30],
Xiong-chaos [40], distributional chaos [23], and so on. These definitions all focus on the complex trajectory behavior
of points. The concept of distributional chaos (it was called strong chaos then) was introduced by Schweizer and Smital
from the perspective of  probability theory for a continuous map in a compact interval [23]. Since then, it has evolved into
three mutually nonequivalent versions of distributional chaos for a map in a metric space: DC1, DC2, and DC3 [3, 28].
DC1 was the original version of distributional chaos introduced in [23], and it is the strongest one among these three
definitions, while DC2 and DC3 are its generalizations.
In our present paper, we pay attention to DC1 for NDSs. So please
remember that distributional chaos always means DC1 in the rest of the present paper without special illustration.
We shall study properties of DC2 and DC3 for NDSs in our forthcoming papers.

For ADSs, many elegant results about distributional chaos have been obtained [1, 14, 16, 17, 20, 21, 34--37, 42].
Most of them are about relations between distributional chaos and other concepts of chaos. For example, Oprocha
showed that neither topologically weak mixing nor Devaney chaos implies distributional chaos using one-sided symbolic
dynamical systems [20]. Then, Li {\it et al.} proved that Devaney chaos implies distributional chaos for continuous maps
with shadowing property in compact metric spaces [16]. Shadowing property (also called pseudo-orbit tracing property),
is a property that a pseudo-orbit can be traced by a true trajectory, which has been extensively studied (see, for example, [6, 15, 39]).
The asymptotic average shadowing property was introduced in [10]. Recently, Wang {\it et al.} proved that the
asymptotic average shadowing property implies distributional chaos for a continuous map with two almost period points
in a compact metric space [37]. Note that Xiong-chaos, introduced in [40], is induced by a topologically mixing map.
In 2002, Yang introduced the notion of distributional chaos in a sequence, and investigated
some relationships of topologically weak mixing, Xiong-chaos and distributional chaos in a sequence for a continuous map
in a locally compact metric space [42]. In 2007, Wang {\it et al.} established a criterion of
distributional chaos under some expanding conditions, and also proved that Li-Yorke chaos is equivalent to distributional
chaos in a sequence for a continuous map in a compact interval [35]. Later, Li and Tan generalized
their result and proved that Li-Yorke $\delta$-chaos is equivalent to distributional $\delta$-chaos
in a sequence for a continuous map in a compact metric space [14].

For NDSs, Kolyada and Snoha extended the concept of topological entropy for ADSs to NDSs and studied its properties [14].
Later, many scholars investigated chaos of NDSs [4, 5, 8, 18, 19, 24, 27, 31--33, 38, 44]. For example, Tian and Chen extended the
concept of Devaney chaos to NDSs in [33]. Then, Chen with the second author in the present paper  generalized related concepts of chaos, such as topological transitivity, sensitivity, chaos in the sense of Li-Yorke, Wiggins, and Devaney to general NDSs, and established a criterion of Li-Yorke chaos induced by strict coupled-expansion for a certain irreducible transitive matrix [27]. In 2012, Balibrea and Oprocha investigated
relations of Li-Yorke chaos with positive topological entropy and topologically weak mixing for NDSs [4]. Recently, we studied whether transitivity and density of periodic points imply sensitivity in the definition of Devaney chaos [44]. However, there are only a few results about distributional chaos for NDSs [8, 18, 38]. To the best of our knowledge, there are no criteria of distributional chaos established for NDSs in the current literature. Motivated by the above works, we shall study distributional chaos in NDSs and try to establish its criteria in the present paper.

The rest of the paper is organized as follows. Section 2 presents some related concepts and useful lemmas.
It is divided into four parts for convenience. In Section 2.1, some basic concepts are given.
In Sections 2.2 and 2.3, several lemmas
about density of a sequence and some properties of one-sided symbolic dynamical systems are recalled,
respectively. They are useful to Sections 3 and 4.  Several kinds of relations for NDSs are
introduced in Section 2.4. Section 3 reveals relations between Li-Yorke $\delta$-chaos and distributional $\delta$-chaos
in a sequence; and Section 4 gives out three criteria of distributional $\delta$-chaos in compact metric spaces, which are caused
by topologically weak mixing, asymptotic average shadowing property, and some expanding condition, respectively. Finally,
a criterion of distributional chaos in a sequence induced by a Xiong chaotic set is established in general metric spaces in Section 5.

\bigskip

\noindent{\bf 2. Preliminaries}\bigskip

In this section, some related concepts and useful lemmas are presented.\medskip

\noindent{\bf 2.1. Some basic concepts}\medskip

For any fixed $x_0\in X$, $\{x_n\}_{n=0}^{\infty}$ is called the (positive) orbit of system (1.1)
starting from $x_0$, where $x_n=f_0^n(x_0)$ and $f_0^n:=f_{n-1}\circ\cdots\circ f_{0}$ for $n\geq 1$.
For convenience, by $f_{0}^{0}$ denote the identity map on $X$, $f_{0}^{-n}:=(f_{0}^{n})^{-1}$,
and $f_{0,\infty}:=\{f_n\}_{n=0}^{\infty}$. By $\overline{A}$ denote the closure of a subset $A\subset X$,
and by $\mathbf{N}$ and $\mathbf{Z^{+}}$ denote the set of all nonnegative integers and that of all positive integers,
respectively.\medskip

Let $P=\{p_n\}_{n=1}^{\infty}\subset\mathbf{N}$ be an increasing sequence, $x,y\in X$, and $\epsilon>0$. Denote
\vspace{-0.2cm}$$F_{x,y}^{*}(\epsilon,P)=
\limsup_{n\to\infty}\frac{1}{n}\sum_{i=1}^{n}\chi_{[0,\epsilon)}\big(d(f_{0}^{p_i}(x), f_{0}^{p_i}(y))\big),\vspace{-0.2cm}$$
\vspace{-0.2cm}$$F_{x,y}(\epsilon,P)=
\liminf_{n\to\infty}\frac{1}{n}\sum_{i=1}^{n}\chi_{[0,\epsilon)}\big(d(f_{0}^{p_i}(x), f_{0}^{p_i}(y))\big),\vspace{-0.2cm}$$
where $\chi_{[0,\epsilon)}$ is the characteristic function defined on the set $[0, \epsilon)$.
$F_{x,y}^{*}$ and $F_{x,y}$ are called the upper and lower probability distributional functions, respectively.\medskip

\noindent{\bf Definition 2.1.} Let $D\subset X$ contain at least two distinct points. Then
$D$ is called a distributionally scrambled set in a sequence $P$ of system (1.1) if,
for any $x, y\in D$ with $x\neq y$,
\begin{itemize}\vspace{-0.2cm}
\item[{\rm (i)}] $F_{x,y}^{*}(\epsilon,P)=1$ for any $\epsilon>0$,\vspace{-0.2cm}
\item[{\rm (ii)}] $F_{x,y}(\delta_{x,y},P)=0$ for some $\delta_{x,y}>0$,
\end{itemize}\vspace{-0.2cm}
and then the pair $(x,y)$ is called a distributionally $\delta_{x,y}$-scrambled pair in the sequence $P$.
Further, $D$ is called a distributionally $\delta$-scrambled set in the sequence $P$
if there exists $\delta>0$ such that (i) holds and $F_{x,y}(\delta,P)=0$ for any $x, y\in D$ with $x\neq y$.
Moreover, if $P=\mathbf{N}$, then $D$ is called a distributionally scrambled set (distributionally $\delta$-scrambled set),
and $F_{x,y}^{*}(\epsilon,\mathbf{N})$ and $F_{x,y}(\delta,\mathbf{N})$ are briefly denoted
by $F_{x,y}^{*}(\epsilon)$ and $F_{x,y}(\delta)$, respectively.\medskip

\noindent{\bf Definition 2.2.} If system (1.1) has an uncountable distributionally scrambled set
(distributionally $\delta$-scrambled set) in a sequence $P$, then it is said to be distributionally chaotic
(distributionally $\delta$-chaotic) in the sequence $P$. Further, if $P=\mathbf{N}$,
then system (1.1) is said to be distributionally chaotic (distributionally $\delta$-chaotic).\medskip

\noindent{\bf Definition 2.3} [27, Definition 2.7]. Let $S\subset X$ contain at least two distinct points.
Then, $S$ is called a Li-Yorke scrambled set of system (1.1) if, for any two distinct points $x, y\in S$,
their corresponding orbits satisfy
\begin{itemize}\vspace{-0.2cm}
\item[{\rm (i)}] $\liminf\limits_{n\rightarrow\infty}d(f_{0}^{n}(x), f_{0}^{n}(y))=0,$ \vspace{-0.2cm}
\item[{\rm (ii)}] $\limsup\limits_{n\rightarrow\infty}d(f_{0}^{n}(x), f_{0}^{n}(y))>0.$ \vspace{-0.2cm}
\end{itemize}
Further, $S$ is called a Li-Yorke $\delta$-scrambled set for some positive constant $\delta$ if, for any two distinct
points $x, y\in S$, (i) holds and, instead of (ii), the following holds:
\begin{itemize}\vspace{-0.2cm}
\item[{\rm (iii)}] $\limsup\limits_{n\rightarrow\infty}d(f_{0}^{n}(x), f_{0}^{n}(y))>\delta.$\vspace{-0.2cm}
\end{itemize}\medskip

\noindent{\bf Definition 2.4}. System (1.1) is said to be Li-Yorke chaotic ($\delta$-chaotic)
if it has an uncountable Li-Yorke scrambled set ($\delta$-scrambled set). Further, system (1.1) is said to be
densely Li-Yorke chaotic ($\delta$-chaotic) if it has a densely uncountable Li-Yorke scrambled set
($\delta$-scrambled set).\medskip

\noindent{\bf Remark 2.5.} Definitions 2.1 and 2.2 generalize the concept of distributional chaos in a sequence
for ADSs introduced in [35, 42] to NDSs. In addition, the concept of Li-Yorke chaos for NDSs in Definition 2.4
was introduced in [27, Definition 2.8].\medskip

\noindent{\bf Definition 2.6} [27, Definition 2.2]. System (1.1) is said to be topologically
transitive in $X$ if for any two nonempty open subsets $U,V\subset X$, there exists $n>0$ such that
$f_0^{n}(U)\cap V\neq\emptyset$.\medskip

It can be easily verified that system {\rm(1.1)} is topologically transitive in $X$ if and only if for any
two nonempty open subsets $U,V\subset X$, there exists $n>0$ such that $f_0^{-n}(U)\cap V\neq\emptyset$.\medskip

\noindent{\bf Definition 2.7.} System (1.1) is said to be topologically weakly mixing in $X$
if for any four nonempty open subsets $U_1, V_1, U_2, V_2$ of $X$, there exists $n>0$ such that
$f_0^{n}(U_1)\cap V_1\neq\emptyset$ and $f_0^{n}(U_2)\cap V_2\neq\emptyset$.\medskip

Note that topologically weak mixing implies topological transitivity for system (1.1).
Next, the concepts of Xiong-chaotic set (see [40]) and asymptotic average shadowing property (see [10])
for ADSs are extended to NDSs.\medskip

\noindent{\bf Definition 2.8.} A subset $S\subset X$ is called a Xiong-chaotic set with respect to a sequence
$P=\{p_n\}_{n=1}^{\infty}\subset\mathbf{N}$ of system (1.1) if for any continuous map $F: S\to X$, there is a subsequence
$\{p_{n_k}\}_{k=1}^{\infty}$ of $P$ such that
\vspace{-0.2cm}$$\lim_{k\to\infty}f_{0}^{p_{n_k}}(x) =F(x),\;\forall\, x\in S.$$

\noindent{\bf Definition 2.9.} A sequence $\{x_n\}_{n=0}^{\infty}\subset X$ is called an asymptotic average pseudo-orbit of system (1.1) if
\vspace{-0.2cm}$$\lim_{n\to\infty}\frac{1}{n}\sum_{i=0}^{n-1}d(f_i(x_i),x_{i+1})=0.\vspace{-0.1cm}$$
The sequence $\{x_n\}_{n=0}^{\infty}$ is asymptotically shadowed in average by a point $y_0\in X$ if
\vspace{-0.2cm}$$\lim_{n\to\infty}\frac{1}{n}\sum_{i=0}^{n-1}d(f_{0}^{i}(y_0),x_{i})=0.$$

\noindent{\bf Definition 2.10.} System (1.1) has asymptotic average shadowing property if any asymptotic average pseudo-orbit of system (1.1)
is asymptotically shadowed in average by some point in $X$.\bigskip

\noindent{\bf 2.2. Density of a sequence}\medskip

In this subsection, some concepts and useful lemmas about density of a sequence are given.\medskip

Let $Q=\{q_n\}_{n=1}^{\infty}\subset\mathbf{N}$ be an increasing sequence
and $P\subset Q$. By $|P|$ denote the cardinality of $P$. Then
$\overline{d}(P|Q):=\limsup\limits_{n\to\infty}\frac{1}{n}|P\cap\{q_1,\dots,q_n\}|$ and
$\underline{d}(P|Q):=\liminf\limits_{n\to\infty}\frac{1}{n}|P\cap\{q_1,\dots,q_n\}|$
are called upper and lower densities of $P$ with respect to $Q$, respectively.
If $\overline{d}(P|Q)=\underline{d}(P|Q)=: d(P|Q)$, then $d(P|Q)$ is called the density of $P$ with respect to $Q$.
If $Q=\mathbf{N}$, then $\overline{d}(P|Q)$, $\underline{d}(P|Q)$,
and $d(P|Q)$ (briefly denoted by $\overline{d}(P)$, $\underline{d}(P)$, and $d(P)$ in this case) are called the upper density,
lower density, and density of $P$, respectively.\medskip

\noindent{\bf Lemma 2.11} [9, Lemma 2]. {\it Let $\{S_n\}_{n=1}^{\infty}$ be a sequence of increasing sequences in $\mathbf{Z^{+}}$.
Then there exists an increasing sequence $Q\subset\mathbf{Z^{+}}$ such that $\overline{d}\bigl((S_n\cap Q)| Q\bigr)=1$ for each $n\geq 1$.}\medskip

\noindent{\bf Lemma 2.12} [22, Lemma 2.6.2]. {\it Let $\{a_n\}_{n=0}^{\infty}$ be a bounded sequence of nonnegative numbers.
Then $\lim\limits_{\scriptstyle n\to\infty}\frac{1}{n}\sum_{i=0}^{n-1}a_i$ $=0$ if and only if there exists a subset
$E\subset \mathbf{N}$ of density zero such that
$\lim\limits_{\scriptstyle i\to\infty \atop \scriptstyle i\notin E}a_i=0$.}
\bigskip

\noindent{\bf 2.3. Some properties of one-sided symbolic dynamical systems}\medskip

In this subsection, we shall briefly recall some properties of one-sided symbolic dynamical systems,
which will be used in subsequent sections. For more details, see [43].\medskip

Let $S=\{0, 1\}$. The one-sided sequence space
$\Sigma_{2}^{+}:=\{\alpha=(a_{0}, a_{1}, a_{2},\cdots):\; a_{i}\in S, \;i\geq 0\}$
is a metric space with the distance
$\rho(\alpha, \beta)=\sum_{i=0}^{\infty}d(a_{i}, b_{i})/2^{i}$,
where $\alpha=(a_{0}, a_{1}, a_{2},\cdots)$, $\beta=(b_{0}, b_{1}, b_{2},\cdots)\in \Sigma_{2}^{+}$,
$d(a_{i}, b_{i})=1$ if $a_{i}\neq b_{i}$, and $d(a_{i}, b_{i})=0$ if $a_{i}=b_{i}$ for $i\geq 0$.
It is a compact metric space. Define the shift map $\sigma: \Sigma_{2}^{+}\rightarrow \Sigma_{2}^{+}$ by
$\sigma(\alpha):=(a_{1}, a_{2},\cdots)$, where $\alpha=(a_{0}, a_{1}, a_{2},\cdots)$. This map is continuous and $(\Sigma_{2}^{+}, \sigma)$ is called the one-sided symbolic dynamical system on two symbols.\medskip

\noindent{\bf Lemma 2.13} [34, Lemma 2.2]. {\it $\Sigma_{2}^{+}$ has an uncountable subset $E$ such that for any different points $\alpha=(a_0,a_1,\cdots),\;\beta=(b_0,b_1,\cdots)$ in $E$, $a_n=b_n$ for infinitely $n$ and $a_m\neq b_m$ for infinitely $m$.}
\bigskip

\noindent{\bf 2.4. Some relations for NDSs}\medskip

In this subsection, several kinds of relations for NDSs are introduced, and
some useful lemmas are presented.\medskip

Let $P=\{p_n\}_{n=1}^{\infty}\subset\mathbf{N}$ be an increasing sequence. The sets
$PR(f_{0,\infty},P):=\{(x,y)\in X\times X: \liminf\limits_{i\to\infty}d(f_{0}^{p_i}(x),f_{0}^{p_i}(y))
=0\}$, $AR(f_{0,\infty},P):=\{(x,y)\in X\times X: \lim\limits_{i\to\infty}d(f_{0}^{p_i}(x),f_{0}^{p_i}$
$(y))=0\}$, $DR(f_{0,\infty},P):=\{(x,y)\in X\times X: \liminf\limits_{i\to\infty}d(f_{0}^{p_i}(x),f_{0}^{p_i}(y))>0\}$
are called the proximal, asymptotic, and distal relations of system (1.1) with respect to $P$, respectively;
the sets $DSR(f_{0,\infty},P):=\{(x,y)\in X\times X: (x,y)$ is a distributionally scrambled pair in $P\}$
and $DSR_{\delta}(f_{0,\infty},P):=\{(x,y)\in X\times X: (x,y)$ is a distributionally $\delta$-scrambled
pair in $P\}$ are called the distributionally chaotic and $\delta$-chaotic relations
of system (1.1) with respect to $P$, respectively; and the sets
$LYR(f_{0,\infty}):=\{(x,y)\in X\times X: (x,y)$ is a Li-Yorke scrambled pair$\}$ and
$LYR_{\delta}(f_{0,\infty}):=\{(x,y)\in X\times X: (x,y)$ is a Li-Yorke $\delta$-scrambled pair $\}$
are called the Li-Yorke chaotic and $\delta$-chaotic relations of system (1.1), respectively.\medskip

\noindent{\bf Proposition 2.14.} {\it Let $P\subset\mathbf{N}$ be an increasing sequence. Then
\begin{itemize}\vspace{-0.2cm}
\item[{\rm (i)}] $DR(f_{0,\infty},P)=(X\times X)\setminus PR(f_{0,\infty},P)$;
\item[{\rm (ii)}] $LYR(f_{0,\infty})=PR(f_{0,\infty},\mathbf{N})\setminus AR(f_{0,\infty},\mathbf{N})$;
\item[{\rm (iii)}] $DSR(f_{0,\infty},P)\subset PR(f_{0,\infty},P)\setminus AR(f_{0,\infty},P)\subset LYR(f_{0,\infty})$.
\end{itemize}}

\noindent{\bf Proof.} It suffices to show Assertion (iii) since Assertions (i) and (ii)
can be directly derived from the above definitions. It is evident that
$PR(f_{0,\infty},P)\setminus AR(f_{0,\infty},P)\subset LYR(f_{0,\infty})$. We shall show that
\vspace{-0.2cm}$$DSR(f_{0,\infty},P)\subset PR(f_{0,\infty},P)\setminus AR(f_{0,\infty},P).                                          \eqno(2.1)\vspace{-0.2cm}$$
Let $P=\{p_n\}_{n=1}^{\infty}$. Fix any $(x,y)\in DSR(f_{0,\infty},P)$.
Then, we have that
\vspace{-0.2cm}$$\limsup_{n\to\infty}\frac{1}{n}\sum_{i=1}^{n}\chi_{[0,\epsilon)}
\bigl(d(f_{0}^{p_i}(x),f_{0}^{p_i}(y))\bigr)=1,\;\forall\;\epsilon>0,               \eqno(2.2)\vspace{-0.2cm}$$
and
\vspace{-0.2cm}$$\liminf_{n\to\infty}\frac{1}{n}\sum_{i=1}^{n}\chi_{[0,\delta)}
\bigl(d(f_{0}^{p_i}(x),f_{0}^{p_i}(y))\bigr)=0                                    \eqno(2.3)\vspace{-0.2cm}$$
for some $\delta>0$. Firstly, we shall show that
\vspace{-0.2cm}$$\liminf_{n\to\infty}d(f_{0}^{p_n}(x),f_{0}^{p_n}(y))=0.          \eqno(2.4)\vspace{-0.2cm}$$
Otherwise, there exists $\xi>0$ such that $\liminf_{n\to\infty}d(f_{0}^{p_n}(x),f_{0}^{p_n}(y))=\xi$.
Then there exists an integer $N>0$ such that $d(f_{0}^{p_n}(x),f_{0}^{p_n}(y))\geq\xi/2$ for each $n\geq N$. Thus,
\vspace{-0.2cm}$$\lim_{n\to\infty}\frac{1}{n}\sum_{i=1}^{n}\chi_{[0,\frac{\xi}{2})}
\bigl(d(f_{0}^{p_i}(x),f_{0}^{p_i}(y))\bigr)=0,\vspace{-0.2cm}$$
which is a contradiction to (2.2).
Secondly, we shall show that
\vspace{-0.2cm}$$\limsup_{n\to\infty}d(f_{0}^{p_n}(x),f_{0}^{p_n}(y))>0.                                                             \eqno(2.5)\vspace{-0.2cm}$$
Otherwise, $\limsup_{n\to\infty}d(f_{0}^{p_n}(x),f_{0}^{p_n}(y))=0$.
So there exists $K>0$ such that $d(f_{0}^{p_n}(x),$ $f_{0}^{p_n}(y))<\delta$ for each $n\geq K$. This implies that
\vspace{-0.2cm}$$\lim_{n\to\infty}\frac{1}{n}\sum_{i=1}^{n}\chi_{[0,\delta)}
\bigl(d(f_{0}^{p_i}(x),f_{0}^{p_i}(y))\bigr)=1,\vspace{-0.2cm}$$
which is a contradiction to (2.3). It follows from (2.4) and (2.5) that $(x,y)\in PR(f_{0,\infty},P)\setminus AR(f_{0,\infty},P)$,
which implies that (2.1) holds.
This completes the proof.\medskip

By (iii) of Proposition 2.14 one can directly get that if system {\rm (1.1)} is distributionally chaotic in a sequence,
then it is Li-Yorke chaotic.\medskip

The following two lemmas are useful to find distributionally or Li-Yorke scrambled sets in Sections 3 and 4.\medskip

\noindent{\bf Lemma 2.15} [11, Lemma 3.2]. {\it Let $R$ be a relation on a complete separable metric space $X$. If $R$ contains
a dense $G_{\delta}$ subset of $X\times X$, then there is a dense $G_{\delta}$ subset $A$ of $X$ such that for each $x\in A$,
there exists a dense $G_{\delta}$ subset $A_{x}$ of $X$ with
$\{(x,y): x\in A, y\in A_{x}\}\subset R$.}\medskip

\noindent{\bf Lemma 2.16} [11, Lemma 3.1]. {\it Assume that $X$ is a complete separable metric space without isolated points. If $R$ is
a symmetric relation with the property that there is a dense $G_{\delta}$ subset $A$ of $X$ such that for each $x\in A$, $R(x)$
contains a dense $G_{\delta}$ subset, then there is a dense subset $B$ of $X$ with uncountably many points such that $(B\times B)\setminus\triangle\subset R$, where $\triangle:=\{(x,x): x\in X\}$.}\medskip

In the above,  $R(x)=\{y: (x,y)\in R\}$ for a relation $R\subset X\times X$ and $x\in X$.

\bigskip

\noindent{\bf 3. Relations between Li-Yorke $\delta$-chaos and distributional $\delta$-chaos in a sequence}\bigskip

In this section, we shall investigate relations between Li-Yorke $\delta$-chaos and distributional $\delta$-chaos in a sequence
for system (1.1).\medskip

\noindent{\bf Lemma 3.1.} {\it Let $(X,d)$ be a separable metric space and $S\subset X$ be an uncountable set.
Then $S$ has an uncountable subset which contains no isolated points.}\medskip

\noindent{\bf Proof.} By the assumptions, $(S,d)$ is a separable metric subspace of $(X,d)$.
Let $\{V_n\}_{n=1}^{\infty}$ be a countable topology basis of $S$. And let $W$ be
the union of all those at most countable sets of $V_n$. Then $W$ is countable, and consequently
 $\tilde{S}:=S\setminus W$ is uncountable. We claim that $\tilde{S}$ contains
 no isolated points. By $\{V_{n_k}\}$ denote the collection of all those uncountable sets of $V_n$.
Then $\{V_{n_k}\cap \tilde{S}\}$ is a topology basis of $\tilde{S}$. Thus, for any $x\in \tilde{S}$
and any neighborhood $U(x)$ of $x$ in $\tilde{S}$, there exists $k_0\geq1$ such that $V_{n_{k_0}}\cap
\tilde{S}\subset U(x)$. Note that $V_{n_{k_0}}\cap \tilde{S}=V_{n_{k_0}}\setminus W$, which implies that
$V_{n_{k_0}}\cap \tilde{S}$ is uncountable. Hence, $x$ is not an isolated
point in $\tilde{S}$, and consequently, $\tilde{S}$ contains no isolated points.
This completes the proof.\medskip

Now, we shall investigate some properties of the set of all distributionally $\delta$-scrambled pairs in a sequence.
The metric of $X\times X$ is denoted by $d'$, and defined by $d'((x_1,x_2), (y_1,y_2))$ $:=\max\{d(x_1,y_1),d(x_2,y_2)\}$
for any $(x_1,x_2), (y_1,y_2)\in X\times X$. Let $W\subset X\times X$ be a nonempty set,
$(x,y)\in X\times X$, and $Q\subset\mathbf{N}$ be an increasing sequence.
Now, we introduce the following notations for convenience:
\vspace{-0.2cm}$$\mathcal{M}_{Q}(1):=\{P\subset Q: \overline{d}(P|Q)=1\},\;N((x,y),W,Q):=\{n\in Q: (f_{0}^{n}(x),f_{0}^{n}(y))\in W\},\vspace{-0.2cm}$$
and
\vspace{-0.2cm}$$\mathcal{F}(W,Q,\mathcal{M}_{Q}(1)):=\{(x,y)\in X\times X: N((x,y),W,Q)\in \mathcal{M}_{Q}(1)\}.$$

It can be easily verified that if $W_1$ and $W_2$ are two nonempty subsets of $X\times X$ with $W_1\subset W_2$, then
\vspace{-0.2cm}$$\mathcal{F}(W_1,Q,\mathcal{M}_{Q}(1))\subset\mathcal{F}(W_2,Q,\mathcal{M}_{Q}(1)).                                  \eqno(3.1)\vspace{-0.2cm}$$

\noindent{\bf Proposition 3.2.} {\it Let $(X,d)$ be a compact metric space and $Q\subset\mathbf{N}$ be an increasing sequence.
Then, $(x,y)\in X\times X$ is a distributionally $\delta$-scrambled pair in $Q$ for some $\delta>0$ if and only if
\begin{itemize}\vspace{-0.2cm}
\item[{\rm (i)}] $(x,y)\in \mathcal{F}([\triangle]_{\epsilon},Q,\mathcal{M}_{Q}(1))$ for all $\epsilon>0$;
\item[{\rm (ii)}] $(x,y)\in \mathcal{F}((X\times X)\setminus\overline{[\triangle]}_{\delta'},Q,\mathcal{M}_{Q}(1))$ for some $\delta'>0$,
\end{itemize}\vspace{-0.2cm}
where $[\triangle]_{\epsilon}=\{(x,y)\in X\times X: d'((x,y),\triangle)<\epsilon\}$.}\medskip

\noindent{\bf Proof.} Let $Q=\{q_n\}_{n=1}^{\infty}$. The proof is divided into two parts.

{\bf(1)  Sufficiency.} Suppose that $(x,y)\in X\times X$ satisfies conditions (i) and (ii). Fix any $\epsilon>0$.
Since $(x,y)\in \mathcal{F}([\triangle]_{\epsilon},Q,\mathcal{M}_{Q}(1))$, one has that
\vspace{-0.2cm}$$\limsup_{n\to\infty}\frac{1}{n}\big|N((x,y),[\triangle]_{\epsilon},Q)\cap\{q_1,\dots,q_n\}\big|=1.                  \eqno(3.2)\vspace{-0.2cm}$$
Let $n\in \mathbf{Z^{+}}$ and $k\in N((x,y),[\triangle]_{\epsilon},Q)\cap\{q_1,\dots,q_n\}$. Then one gets that
\vspace{-0.2cm}$$d'((f_{0}^{k}(x),f_{0}^{k}(y)),\triangle)<\epsilon.\vspace{-0.2cm}$$
Noting that $\triangle$ is a closed subset in compact metric space $(X\times X,d')$,
we have that there exists $(z,z)\in \triangle$ such that
\vspace{-0.2cm}$$d'((f_{0}^{k}(x),f_{0}^{k}(y)),\triangle)=d'((f_{0}^{k}(x),f_{0}^{k}(y)),(z,z))<\epsilon,  \vspace{-0.2cm}$$
and thus
\vspace{-0.2cm}$$d(f_{0}^{k}(x),f_{0}^{k}(y))<2\epsilon.\vspace{-0.2cm}$$
So by (3.2) one has that
\vspace{-0.2cm}$$\limsup_{n\to\infty}\frac{1}{n}\big|\big\{k\in\{q_1,\dots,q_n\}: d(f_{0}^{k}(x),f_{0}^{k}(y))<2\epsilon\big\}\big|=1,\vspace{-0.2cm}$$
which implies that
\vspace{-0.2cm}$$F_{x,y}^{*}(2\epsilon,Q)=1.                                                                                         \eqno(3.3)\vspace{-0.2cm}$$
On the other hand, since $(x,y)\in \mathcal{F}((X\times X)\setminus\overline{[\triangle]}_{\delta'},Q,\mathcal{M}_{Q}(1))$, we have that
\vspace{-0.2cm}$$\limsup_{n\to\infty}\frac{1}{n}\big|N((x,y),(X\times X)\setminus\overline{[\triangle]}_{\delta'},Q)\cap\{q_1,\dots,q_n\}\big|=1.                   \eqno(3.4)\vspace{-0.2cm}$$
Let $n\in \mathbf{Z^{+}}$ and $k\in N((x,y),(X\times X)\setminus\overline{[\triangle]}_{\delta'},Q)\cap\{q_1,\dots,q_n\}$. Then we obtain that
\vspace{-0.2cm}$$d'((f_{0}^{k}(x),f_{0}^{k}(y)),\triangle)>\delta',\vspace{-0.2cm}$$
and thus
\vspace{-0.2cm}$$d(f_{0}^{k}(x),f_{0}^{k}(y))=d'((f_{0}^{k}(x),f_{0}^{k}(y)),(f_{0}^{k}(x),f_{0}^{k}(x)))>\delta'.\vspace{-0.2cm}$$
This, together with (3.4), implies that
\vspace{-0.2cm}$$\limsup_{n\to\infty}\frac{1}{n}\big|\big\{k\in\{q_1,\dots,q_n\}: d(f_{0}^{k}(x),f_{0}^{k}(y))>\delta'\big\}\big|=1.\vspace{-0.2cm}$$
Hence,
\vspace{-0.2cm}$$\liminf_{n\to\infty}\frac{1}{n}\big|\big\{k\in\{q_1,\dots,q_n\}: d(f_{0}^{k}(x),f_{0}^{k}(y))<\delta'\big\}\big|=0,\vspace{-0.2cm}$$
which results in
\vspace{-0.2cm}$$F_{x,y}(\delta',Q)=0.                                                                                               \eqno(3.5)\vspace{-0.2cm}$$
Therefore, it follows from (3.3) and (3.5) that $(x,y)$ is a distributionally $\delta'$-scrambled pair in the sequence $Q$.

{\bf(2)  Necessity.} Suppose that $(x,y)\in X\times X$ is a distributionally $\delta$-scrambled pair in the sequence $Q$ for some $\delta>0$.
Then, for any $\epsilon>0$, one has that
\vspace{-0.2cm}$$\limsup_{n\to\infty}\frac{1}{n}\big|\big\{k\in\{q_1,\dots,q_n\}: d(f_{0}^{k}(x),f_{0}^{k}(y))<\epsilon\big\}\big|=1;        \eqno(3.6)\vspace{-0.2cm}$$
and
\vspace{-0.2cm}$$\liminf_{n\to\infty}\frac{1}{n}\big|\big\{k\in\{q_1,\dots,q_n\}: d(f_{0}^{k}(x),f_{0}^{k}(y))<\delta\big\}\big|=0.          \eqno(3.7)\vspace{-0.2cm}$$
Let $n\in \mathbf{Z^{+}}$ and $k\in\{q_1,\dots,q_n\}$ such that $d(f_{0}^{k}(x),f_{0}^{k}(y))<\epsilon$.
Then
\vspace{-0.2cm}$$d'((f_{0}^{k}(x),f_{0}^{k}(y)),(f_{0}^{k}(x),f_{0}^{k}(x)))=d((f_{0}^{k}(x),f_{0}^{k}(y)))<\epsilon,\vspace{-0.2cm}$$
which implies that $(f_{0}^{k}(x),f_{0}^{k}(y))\in [\triangle]_{\epsilon}$.
Thus, by (3.6) we get that
\vspace{-0.2cm}$$\limsup_{n\to\infty}\frac{1}{n}\big|\big\{k\in\{q_1,\dots,q_n\}: (f_{0}^{k}(x),f_{0}^{k}(y))\in [\triangle]_{\epsilon}\big\}\big|=1.\vspace{-0.2cm}$$
So $(x,y)\in \mathcal{F}([\triangle]_{\epsilon},Q,\mathcal{M}_{Q}(1))$ for any $\epsilon>0$.
On the other hand, it follows from (3.7) that
\vspace{-0.2cm}$$\limsup_{n\to\infty}\frac{1}{n}\big|\big\{k\in\{q_1,\dots,q_n\}: d(f_{0}^{k}(x),f_{0}^{k}(y))\geq\delta\big\}\big|=1.             \eqno(3.8)\vspace{-0.2cm}$$
Let $n\in \mathbf{Z^{+}}$ and $k\in\{q_1,\dots,q_n\}$ such that $d(f_{0}^{k}(x),f_{0}^{k}(y))\geq\delta$.
We claim that $(f_{0}^{k}(x),f_{0}^{k}(y))$ $\in (X\times X)\setminus\overline{[\triangle]}_{\delta/4}$.
Otherwise, $d'((f_{0}^{k}(x),f_{0}^{k}(y)),\triangle)\leq\delta/4$. Then there exists $(z,z)\in \triangle$ such that
\vspace{-0.2cm}$$d'((f_{0}^{k}(x),f_{0}^{k}(y)),\triangle)=d'((f_{0}^{k}(x),f_{0}^{k}(y)),(z,z))\leq\delta/4.\vspace{-0.2cm}$$
So $d(f_{0}^{k}(x),f_{0}^{k}(y))\leq\delta/2$, which is a contradiction. Thus, by (3.8) one gets that
\vspace{-0.2cm}$$\limsup_{n\to\infty}\frac{1}{n}\big|\big\{k\in\{q_1,\dots,q_n\}: (f_{0}^{k}(x),f_{0}^{k}(y))\in
(X\times X)\setminus\overline{[\triangle]}_{\delta/4}\big\}\big|=1.               \vspace{-0.2cm}$$
Then $(x,y)\in \mathcal{F}((X\times X)\setminus\overline{[\triangle]}_{\delta'},Q,\mathcal{M}_{Q}(1))$,
where $\delta':=\delta/4$. Therefore, the entire proof is complete.\medskip

The following result can be directly derived from Proposition 3.2.\medskip

\noindent{\bf Corollary 3.3.} {\it Let $(X,d)$ be a compact metric space, $Q\subset\mathbf{N}$ be an increasing sequence,
and $\delta>0$. Then there exists $\delta'>0$ such that
\vspace{-0.2cm}$$\left(\bigcap_{\epsilon>0}\mathcal{F}([\bigtriangleup]_{\epsilon},Q,\mathcal{M}_{Q}(1))\right)
\bigcap\mathcal{F}((X\times X)\setminus\overline{[\triangle]}_{\delta},Q,\mathcal{M}_{Q}(1))
\subset DSR_{\delta'}(f_{0,\infty},Q).\vspace{-0.2cm}$$}

\noindent{\bf Proposition 3.4.} {\it Let $Q\subset\mathbf{N}$ be an increasing sequence and $f_n: X\to X$
be continuous for each $n\geq0$. Then $\mathcal{F}(W,Q,\mathcal{M}_{Q}(1))$ is a $G_\delta$ subset
of $X\times X$ for any nonempty open subset $W$ of $X\times X$.}\medskip

\noindent{\bf Proof.} Let $Q=\{q_n\}_{n=1}^{\infty}$. For any $k,n\in \mathbf{Z^{+}}$ and any $m\in \mathbf{Z^{+}}$
with $m\geq n+1$, denote
\vspace{-0.2cm}$$\Theta_{k,n,m}:=\{1\leq l\leq m: \frac{l}{m}>1-\frac{1}{k}\}.\vspace{-0.2cm}$$
For any $l\in\Theta_{k,n,m}$, set
\vspace{-0.2cm}$$\Pi_{k,n,m,l}:=\big\{(r_1,\cdots,r_l): \{r_1,\cdots,r_l\}\subset\{q_1,\cdots, q_m\}\big\}.\vspace{-0.2cm}$$
Fix any nonempty open subset $W\subset X\times X$. Denote
\vspace{-0.2cm}$$V:=\bigcap_{k=1}^{\infty}\bigcap_{n=1}^{\infty}\bigcup_{m=n+1}^{\infty}
\bigcup_{l\in\Theta_{k,n,m}}\bigcup_{(r_1,\cdots,r_l)\in\Pi_{k,n,m,l}}\bigcap_{i=1}^{l}(f_{0}^{-r_i}\times f_{0}^{-r_i})(W).                                                         \eqno(3.9)\vspace{-0.2cm}$$

Now, it is to show that $\mathcal{F}(W,Q,\mathcal{M}_{Q}(1))=V$.
On the one hand, for any given $(x,y)\in \mathcal{F}(W,Q,\mathcal{M}_{Q}(1))$, it follows that
\vspace{-0.2cm}$$\limsup\limits_{m\to\infty}\frac{1}{m}\big|\big\{j\in \{q_1,\cdots, q_m\}: (x,y)\in (f_{0}^{-j}\times f_{0}^{-j})(W)\big\}\big|=1.\vspace{-0.2cm}$$
Then, for each $k\in \mathbf{Z^{+}}$, we have that
\vspace{-0.2cm}$$\limsup\limits_{m\to\infty}\frac{1}{m}\big|\big\{j\in \{q_1,\cdots, q_m\}: (x,y)\in (f_{0}^{-j}\times f_{0}^{-j})(W)\big\}\big|>1-\frac{1}{k}.\vspace{-0.2cm}$$
So, for each $n\in\mathbf{Z^{+}}$, there exists $m\geq n+1$ such that
\vspace{-0.2cm}$$\frac{1}{m}\big|\big\{j\in \{q_1,\cdots, q_m\}: (x,y)\in (f_{0}^{-j}\times f_{0}^{-j})(W)\big\}\big|>1-\frac{1}{k}.\vspace{-0.2cm}$$
Thus, there exist $1\leq l\leq m$ with
\vspace{-0.2cm}$$\frac{l}{m}>1-\frac{1}{k},\vspace{-0.2cm}$$
and $(r_1,\cdots,r_l)$ with $\{r_1,\cdots,r_l\}\subset\{q_1,\cdots, q_m\}$ such that
\vspace{-0.2cm}$$(x,y)\in(f_{0}^{-r_i}\times f_{0}^{-r_i})(W),\;1\leq i\leq l.\vspace{-0.2cm}$$
Hence, $(x,y)\in V$, and consequently, $\mathcal{F}(W,Q,\mathcal{M}_{Q}(1))\subset V$.

On the other hand, for any given $(x,y)\in V$, and any $k, n\in\mathbf{Z^{+}},$
there exist $m_n\geq n+1$ and $1\leq l_n\leq m_n$ with
\vspace{-0.2cm}$$\frac{l_n}{m_n}>1-\frac{1}{k},\vspace{-0.2cm}$$
and $(r_1,\cdots,r_{l_n})$ with $\{r_1,\cdots,r_{l_n}\}\subset\{q_1,\cdots, q_{m_n}\}$ such that
\vspace{-0.2cm}$$(x,y)\in(f_{0}^{-r_i}\times f_{0}^{-r_i})(W),\,\;1\leq i\leq l_n,\vspace{-0.2cm}$$
which implies that
\vspace{-0.2cm}$$\frac{1}{m_n}\big|\big\{j\in \{q_1,\cdots, q_{m_n}\}: (x,y)\in (f_{0}^{-j}\times f_{0}^{-j})(W)\big\}\big|\geq\frac{l_n}{m_n}>1-\frac{1}{k}.\vspace{-0.2cm}$$
Then
\vspace{-0.2cm}$$\limsup\limits_{n\to\infty}\frac{1}{m_n}\big|\big\{j\in \{q_1,\cdots, q_{m_n}\}: (x,y)\in (f_{0}^{-j}\times f_{0}^{-j})(W)\big\}\big|\geq1-\frac{1}{k}.\vspace{-0.2cm}$$
Since $k$ is arbitrary, we obtain that
\vspace{-0.2cm}$$\limsup\limits_{n\to\infty}\frac{1}{m_n}\big|\big\{j\in \{q_1,\cdots, q_{m_n}\}: (x,y)\in (f_{0}^{-j}\times f_{0}^{-j})(W)\big\}\big|=1.\vspace{-0.2cm}$$
Hence, $(x,y)\in\mathcal{F}(W,Q,\mathcal{M}_{Q}(1))$, and consequently, $V\subset\mathcal{F}(W,Q,\mathcal{M}_{Q}(1))$.

Therefore, $\mathcal{F}(W,Q,\mathcal{M}_{Q}(1))=V$. Since $f_n$ is continuous in $X$ for each $n\geq0$,
$(f_{0}^{-j}\times f_{0}^{-j})(W)$ is an open subset of $X\times X$ for each $j\geq0$. Thus,
$\mathcal{F}(W,Q,\mathcal{M}_{Q}(1))$ is a $G_\delta$ subset of $X\times X$ by (3.9).
This completes the proof.\medskip

\noindent{\bf Remark 3.5.} The proof of Proposition 3.4 is motivated by that of Theorem 3.2 in [41].\medskip

\noindent{\bf Theorem 3.6.} {\it Let $(X,d)$ be a compact metric space and $f_n: X\to X$ be continuous for each $n\geq0$.
Then system {\rm(1.1)} is Li-Yorke $\delta$-chaotic for some $\delta>0$ if and only if it is
distributionally $\delta'$-chaotic in a sequence for some $\delta'>0$.}\medskip

\noindent{\bf Proof.} The proof is divided into two parts. Although the idea used in the proof of the sufficiency part of Theorem 3.6
is similar to that of Assertion (iii) in Proposition 2.14, we shall give its detailed proof for completeness.

{\bf(1)  Sufficiency.} Suppose that system {\rm(1.1)} is distributionally $\delta'$-chaotic in a sequence $P=\{p_i\}_{i=1}^{\infty}$
for some $\delta'>0$. Let $D\subset X$ be an uncountable distributionally $\delta'$-scrambled set in the sequence $P$.
Then, for any $x,y\in D$ with $x\neq y$, we have that
\vspace{-0.2cm}$$\limsup_{n\to\infty}\frac{1}{n}\sum_{i=1}^{n}\chi_{[0,\epsilon)}
\bigl(d(f_{0}^{p_i}(x),f_{0}^{p_i}(y))\bigr)=1,\;\forall\;\epsilon>0,       \eqno(3.10)\vspace{-0.2cm}$$
and
\vspace{-0.2cm}$$\liminf_{n\to\infty}\frac{1}{n}\sum_{i=1}^{n}\chi_{[0,\delta')}
\bigl(d(f_{0}^{p_i}(x),f_{0}^{p_i}(y))\bigr)=0.                               \eqno(3.11)\vspace{-0.2cm}$$

We shall show that $D$ is a Li-Yorke $\delta$-scrambled set for system (1.1) with $\delta=\delta'/2$.

Firstly, we shall show that for any $x,y\in D$ with $x\neq y$,
\vspace{-0.2cm}$$\liminf_{n\to\infty}d(f_{0}^{p_n}(x),f_{0}^{p_n}(y))=0.                                                            \eqno(3.12)\vspace{-0.2cm}$$
Otherwise, there exist $x_0, y_0\in D$ with $x_0\neq y_0$ and $\xi>0$ such that
\vspace{-0.2cm}$$\liminf_{n\to\infty}d(f_{0}^{p_n}(x_0),f_{0}^{p_n}(y_0))=\xi.\vspace{-0.2cm}$$
Then there exists an integer $N>0$ such that
\vspace{-0.2cm}$$d(f_{0}^{p_n}(x_0),f_{0}^{p_n}(y_0))\geq\frac{\xi}{2}, \;n\geq N.\vspace{-0.2cm}$$
Thus,
\vspace{-0.2cm}$$\lim_{n\to\infty}\frac{1}{n}\sum_{i=1}^{n}\chi_{[0,\frac{\xi}{2})}
\bigl(d(f_{0}^{p_i}(x_0),f_{0}^{p_i}(y_0))\bigr)=0,\vspace{-0.2cm}$$
which is a contradiction to (3.10). So by (3.12) one has that for any $x,y\in D$ with $x\neq y$,
\vspace{-0.2cm}$$\liminf_{n\to\infty}d(f_{0}^{n}(x),f_{0}^{n}(y))=0.                                                                \eqno(3.13)\vspace{-0.2cm}$$

Secondly, we shall show that for any $x,y\in D$ with $x\neq y$,
\vspace{-0.2cm}$$\limsup_{n\to\infty}d(f_{0}^{p_n}(x),f_{0}^{p_n}(y))>\delta.                                                       \eqno(3.14)\vspace{-0.2cm}$$
Otherwise, there exist $x'_0,y'_0\in D$ with $x'_0\neq y'_0$ such that
\vspace{-0.2cm}$$\limsup_{n\to\infty}d(f_{0}^{p_n}(x'_0),f_{0}^{p_n}(y'_0))\leq\delta.\vspace{-0.2cm}$$
So there exists $K>0$ such that
\vspace{-0.2cm}$$d(f_{0}^{p_n}(x'_0),f_{0}^{p_n}(y'_0))<\delta',\;n\geq K.\vspace{-0.2cm}$$
Thus,
\vspace{-0.2cm}$$\lim_{n\to\infty}\frac{1}{n}\sum_{i=1}^{n}\chi_{[0,\delta')}
\bigl(d(f_{0}^{p_i}(x'_0),f_{0}^{p_i}(y'_0))\bigr)=1,$$
which is a contradiction to (3.11). So by (3.14) we get that for any $x,y\in D$ with $x\neq y$,
\vspace{-0.2cm}$$\limsup_{n\to\infty}d(f_{0}^{n}(x),f_{0}^{n}(y))>\delta.       \vspace{-0.2cm}      \eqno(3.15)$$
It follows from (3.13) and (3.15) that $D$ is a Li-Yorke $\delta$-scrambled set.
Therefore, system {\rm(1.1)} is Li-Yorke $\delta$-chaotic.

{\bf(2)  Necessity.} Suppose that system (1.1) is Li-Yorke $\delta$-chaotic and $S$ is an uncountable Li-Yorke
$\delta$-scrambled set for system (1.1). Then $S$ has an uncountable subset $\tilde{S}$ that
has no isolated points by Lemma 3.1. Let $Y:=\overline{\tilde{S}}$. Then, $(Y,d)$ is a compact
metric subspace of $(X,d)$ with no isolated points. We can choose a countable dense subset $C_0$ of $\tilde{S}$
since it is separable. Let $C_0=\{x_i\}_{i=1}^{\infty}$. Then, for any $i\neq j$,
there exist two increasing sequences $P_{ij}=\{t_{n}^{ij}\}_{n=1}^{\infty}\subset\mathbf{Z^{+}}$
and $P'_{ij}=\{m_{n}^{ij}\}_{n=1}^{\infty}\subset\mathbf{Z^{+}}$ such that
\vspace{-0.2cm}$$\lim_{n\to\infty}d(f_{0}^{t_{n}^{ij}}(x_i),f_{0}^{t_{n}^{ij}}(x_j))=0,\;\;
\lim_{n\to\infty}d(f_{0}^{m_{n}^{ij}}(x_i),f_{0}^{m_{n}^{ij}}(x_j))>\delta.                                                         \eqno(3.16)\vspace{-0.2cm}$$
By Lemma 2.11 there exists an increasing sequence $Q\subset\mathbf{Z^{+}}$ such that for any $i\neq j$,
\vspace{-0.2cm}$$\overline{d}(P_{ij}\cap Q|Q)=1,\;\;\overline{d}(P'_{ij}\cap Q|Q)=1.                                                \eqno(3.17)\vspace{-0.2cm}$$
It follows from the first relation in (3.16) that for any $\epsilon>0$, there exists an integer $N_{ij}>0$ such that
\vspace{-0.2cm}$$d(f_{0}^{t_{n}^{ij}}(x_i),f_{0}^{t_{n}^{ij}}(x_j))<\epsilon,\;n>N_{ij},\vspace{-0.2cm}$$
which implies that
\vspace{-0.2cm}$$(f_{0}^{t_{n}^{ij}}(x_i),f_{0}^{t_{n}^{ij}}(x_j))\in [\triangle]_{\epsilon},\;n>N_{ij}.\vspace{-0.2cm}$$
Thus,
\vspace{-0.2cm}$$\{t_{n}^{ij}\}_{n=N_{ij}+1}^{\infty}\cap Q\subset N((x_i,x_j),[\triangle]_{\epsilon},Q),\vspace{-0.2cm}$$
which, together with the first relation in (3.17), yields that
\vspace{-0.2cm}$$(x_i, x_j)\in \mathcal{F}([\bigtriangleup]_{\epsilon},Q,\mathcal{M}_{Q}(1))\; {\rm for}\; {\rm all}\; \epsilon>0.  \eqno(3.18)\vspace{-0.2cm}$$
In addition, it follows from the second relation in (3.16) that there exists an integer $K_{ij}>0$ such that
\vspace{-0.2cm}$$d(f_{0}^{m_{n}^{ij}}(x_i),f_{0}^{m_{n}^{ij}}(x_j))>\delta,\;n>K_{ij},\vspace{-0.2cm}$$
which implies that
\vspace{-0.2cm}$$(f_{0}^{m_{n}^{ij}}(x_i),f_{0}^{m_{n}^{ij}}(x_j))\in (X\times X)\setminus\overline{[\triangle]}_{\frac{\delta}{2}}.\vspace{-0.2cm}$$
Thus,
\vspace{-0.2cm}$$\{m_{n}^{ij}\}_{n=K_{ij}+1}^{\infty}\cap Q\subset N((x_i,x_j),(X\times X)\setminus\overline{[\triangle]}_{\frac{\delta}{2}},Q),\vspace{-0.2cm}$$
which, together with the second relation in (3.17), yields that
\vspace{-0.2cm}$$(x_i, x_j)\in \mathcal{F}((X\times X)\setminus\overline{[\triangle]}_{\frac{\delta}{2}},Q,\mathcal{M}_{Q}(1)).     \eqno(3.19)\vspace{-0.2cm}$$
By (3.18) and (3.19) one has that
\vspace{-0.2cm}$$\begin{array}{llll}(C_0\times C_{0})\setminus\triangle &\subset \left(\bigcap\limits_{\epsilon>0}\mathcal{F}([\bigtriangleup]_{\epsilon},Q,\mathcal{M}_{Q}(1))\right)
\bigcap\mathcal{F}\big((X\times X)\setminus\overline{[\triangle]}_{\frac{\delta}{2}},Q,\mathcal{M}_{Q}(1)\big).
\end{array}\vspace{-0.2cm}$$
Denote
\vspace{-0.2cm}
$$\Omega:=\left(\bigcap_{n=1}^{\infty}\mathcal{F}([\bigtriangleup]_{\frac{1}{n}},Q,\mathcal{M}_{Q}(1))\right)
\bigcap\mathcal{F}\big((X\times X)\setminus\overline{[\triangle]}_{\frac{\delta}{2}},Q,\mathcal{M}_{Q}(1)\big).\vspace{-0.2cm}$$
Then $(C_0\times C_0)\setminus\bigtriangleup\subset \Omega$ and $\Omega$ is a $G_{\delta}$ subset of $X\times X$ by
Proposition 3.4. Thus, $\Omega\cap (Y\times Y)$ is a dense $G_{\delta}$ subset of $Y\times Y$ since $C_0$ is dense in $Y$.

By (3.1) it can be easily verified that
\vspace{-0.2cm}
$$\bigcap_{n=1}^{\infty}\mathcal{F}\big([\bigtriangleup]_{\frac{1}{n}},Q,\mathcal{M}_{Q}(1)\big)=
\bigcap_{\epsilon>0}\mathcal{F}\big([\bigtriangleup]_{\epsilon},Q,\mathcal{M}_{Q}(1)\big).\vspace{-0.2cm}$$
Hence,
\vspace{-0.2cm}$$\Omega=\left(\bigcap_{\epsilon>0}\mathcal{F}([\bigtriangleup]_{\epsilon},Q,\mathcal{M}_{Q}(1))\right)
\bigcap\mathcal{F}\big((X\times X)\setminus\overline{[\triangle]}_{\frac{\delta}{2}},Q,\mathcal{M}_{Q}(1)\big).\vspace{-0.2cm}$$
Consequently, there exists $\delta'>0$ such that $\Omega\subset DSR_{\delta'}(f_{0,\infty},Q)$ by Corollary 3.3.
So $DSR_{\delta'}(f_{0,\infty},Q)$ contains a dense $G_{\delta}$ subset of $Y\times Y$.
By Lemmas 2.15 and 2.16, there is a dense subset $B_0$ of $Y$ with uncountably many points such that
$(B_0\times B_0)\setminus\triangle\subset DSR_{\delta'}(f_{0,\infty},Q)$. Hence,
system (1.1) is distributionally $\delta'$-chaotic in the sequence $Q$.
The whole proof is complete.\medskip

\noindent{\bf Remark 3.7.} The result of Theorem 3.6 extends that of Theorem 3 in [14] for ADSs to NDSs.\medskip

\bigskip

\noindent{\bf 4. Three criteria of distributional chaos in compact metric spaces}\bigskip

In this section, we shall establish three criteria of distributional chaos for system (1.1) in the case that the metric space is compact,
which are caused by topologically weak mixing, asymptotic average shadowing property, and some expanding condition, respectively.\medskip

\noindent{\bf Lemma 4.1} [26, Lemma 2.7]. {\it Let $(X,d)$ be a complete metric space and $\{A_n\}_{n=1}^{\infty}$ be a sequence of
bounded and closed sets of $X$ which have the finite intersection property. If the diameter $d(A_n)\to 0$ as $n\to\infty$, then
$\bigcap_{n=1}^{\infty}A_n$ is a singleton set.}\medskip

It was shown that a topologically weakly mixing NDS contains a scrambled set in [4, Theorem 7]. The following result improves this result.\medskip

\noindent{\bf Proposition 4.2.} {\it Let $(X,d)$ be a complete separable metric space without isolated points
and $f_n: X\to X$ be continuous for each $n\geq0$. If system {\rm(1.1)} is topologically weakly mixing,
then it is densely Li-Yorke $\delta$-chaotic for some $\delta>0$.}\medskip

\noindent{\bf Proof.} Denote
\vspace{-0.2cm}$$S:=\{(x,y)\in X\times X: \{(f_{0}^{n}(x),f_{0}^{n}(y))\}_{n=0}^{\infty} {\rm \;is\; dense\; in\;} X\times X\}.      \eqno(4.1)\vspace{-0.2cm}$$
Since $(X,d)$ is a complete separable metric space, $X\times X$ has a countable topology basis $\{U_n\times V_n\}_{n=1}^{\infty}$.
It can be easily verified that
\vspace{-0.2cm}$$S=\bigcap_{n=1}^{\infty}\bigcup_{k\geq0}(f_{0}^{-k}\times f_{0}^{-k})(U_n\times V_n).                         \vspace{-0.2cm}$$
Because system (1.1) is topologically weakly mixing, $\bigcup_{k\geq0}(f_{0}^{-k}\times f_{0}^{-k})(U_n\times V_n)$ is dense in $X\times X$
for each $n\geq1$. Together with the assumption that $f_n$ is continuous for each $n\geq0$, we get that $S$ is a dense $G_{\delta}$ subset of $X\times X$.

Fix any two different points $x_0, y_0\in X$. Let $\delta:=d(x_0,y_0)/2>0$. By (4.1), for any $(x,y)\in S$, there exist two increasing sequences $\{p_n\}_{n=1}^{\infty},\{q_n\}_{n=1}^{\infty}\subset\mathbf{N}$ such that $(f_{0}^{p_n}(x),f_{0}^{p_n}(y))\to(x_0,x_0)$ and $(f_{0}^{q_n}(x),f_{0}^{q_n}(y))\to(x_0,y_0)$ as $n\to\infty$, respectively, which implies that
\vspace{-0.2cm}$$\liminf_{n\to\infty}d(f_{0}^{n}(x),f_{0}^{n}(y))=0,\;\;
\limsup_{n\to\infty}d(f_{0}^{n}(x),f_{0}^{n}(y))\geq d(x_0,y_0)>\delta.\vspace{-0.2cm}$$
So $S\subset LYR_{\delta}(f_{0,\infty})$, and thus $LYR_{\delta}(f_{0,\infty})$ contains a dense $G_{\delta}$ subset of $X\times X$.
By Lemmas 2.15 and 2.16 there is a dense subset $B'$ of $X$ with uncountably many points such that
$(B'\times B')\setminus\triangle\subset LYR_{\delta}(f_{0,\infty})$.
Therefore, system {\rm(1.1)} is densely Li-Yorke $\delta$-chaotic. This completes the proof.\medskip

The following result is a direct consequence of Theorem 3.6 and Proposition 4.2.\medskip

\noindent{\bf Theorem 4.3.} {\it Let $(X,d)$ be a compact metric space without isolated points
and $f_n: X\to X$ be continuous for each $n\geq0$. If system {\rm(1.1)} is topologically weakly mixing,
then it is distributionally $\delta$-chaotic in a sequence for some $\delta>0$.}\medskip

The next result is a criterion of distributional $\delta$-chaos induced by asymptotic average shadowing property
for system (1.1).\medskip

\noindent{\bf Theorem 4.4.} {\it Let $(X,d)$ be a compact metric space and $f_n: X\to X$ be continuous for each $n\geq0$.
Assume that system {\rm (1.1)} has two periodic points whose orbits do not intersect each other. If system {\rm (1.1)}
has asymptotic average shadowing property, then it is distributionally $\delta$-chaotic for some $\delta>0$.}\medskip

\noindent{\bf Remark 4.5.} It is evident that two different periodic orbits of ADSs do not intersect each other.
However, it is not true for NDSs in general (see Example 2.1 in [44]).\medskip

\noindent{\bf Proof.} For convenience, we shall divide the whole proof into two steps.

{\bf Step $1$}. Construct an uncountable set $D_0$.

Suppose that $x$ and $y$ are two periodic points of system (1.1) whose orbits do not intersect each other. Then there exists $\delta>0$ such that
\vspace{-0.2cm}$$d(f_{0}^{n}(x),f_{0}^{n}(y))>2\delta,\; n\geq0.                                                            \eqno(4.2)\vspace{-0.2cm}$$
Fix any countable set $X_0=\{e_i\}_{i=0}^{\infty}$ of $X$. Denote
\vspace{-0.2cm}$$\mathcal{E}:=\{h=(h_{0},h_{1},\cdots): h_{i}\in \{x,y\},\; i\geq0\}.                                       \eqno(4.3)\vspace{-0.2cm}$$
It is clear that $\mathcal{E}$ is uncountable. Define a map $g: \mathcal{E}\to g(\mathcal{E})$ by
\vspace{-0.2cm}$$g(h):=(e_0,h_{0},e_0,e_1,h_{0},h_{1},e_0,e_1,e_2,h_{0},h_{1},h_2,\cdots),\;
h=(h_{0},h_{1},\cdots)\in \mathcal{E}.                                                \eqno(4.4)\vspace{-0.2cm}$$
Evidently, $g$ is bijective, and thus $g(\mathcal{E})$ is uncountable. Let $m_0:=1$ and
\vspace{-0.2cm}$$m_{n+1}:=(2^{n}+1)m_n,\;\;n\geq0.                                                                                \vspace{-0.2cm}$$
For any $u=(u_0,u_1,\cdots)\in g(\mathcal{E})$, we define a sequence $\{x_i(u)\}_{i=0}^{\infty}$ by
$x_0(u)=u_0$ and
\vspace{-0.2cm}$$x_i(u)=f_{0}^{i}(u_n),\;m_n\leq i< m_{n+1},\;\;n\geq0.                \eqno(4.5)\vspace{-0.2cm}$$
It can be easily verified that $\{x_i(u)\}_{i=0}^{\infty}$ is an asymptotic average pseudo-orbit of system (1.1), which
is asymptotically shadowed in average by some point $y_u$. Denote
\vspace{-0.2cm}$$D_0=\{y_u: u\in g(\mathcal{E})\}.                                                                            \vspace{-0.2cm}$$

Fix any two different points
\vspace{-0.2cm}$$u=(u_0,u_1,\cdots),v=(v_0,v_1,\cdots)\in g(\mathcal{E}).                                                      \eqno(4.6)\vspace{-0.2cm}$$
Then
\vspace{-0.2cm}$$\lim_{n\to\infty}\frac{1}{n}\sum_{i=0}^{n-1}d(f_{0}^{i}(y_u),x_{i}(u))=0,\;\;
\lim_{n\to\infty}\frac{1}{n}\sum_{i=0}^{n-1}d(f_{0}^{i}(y_{v}),x_{i}(v))=0.                                                          \eqno(4.7)\vspace{-0.2cm}$$
We claim that $y_u\neq y_{v}$. In fact, suppose that $y_u=y_{v}$. Then, by (4.7) one has that
\vspace{-0.2cm}$$\lim_{n\to\infty}\frac{1}{n}\sum_{i=0}^{n-1}d(x_{i}(u),x_{i}(v))=0.                                                 \eqno(4.8)\vspace{-0.2cm}$$
Since $u\neq v\in g(\mathcal{E})$, it follows from (4.3) and (4.4) that there exists an increasing sequence
$\{l_n\}_{n=0}^{\infty}\subset\mathbf{Z^{+}}$
such that $u_{l_n}=x, v_{l_n}=y$ or $u_{l_n}=y, v_{l_n}=x$ for each $n\geq0$.
Then, by (4.2) and (4.5) we get that
\vspace{-0.2cm}$$d(x_i(u),x_i(v))=d(f_{0}^{i}(u_{l_n}),f_{0}^{i}(v_{l_n}))=d(f_{0}^{i}(x),f_{0}^{i}(y))>2\delta,\;\;m_{l_n}\leq i< m_{l_{n}+1},                                                                                                      \eqno(4.9)\vspace{-0.2cm}$$
which implies that
\vspace{-0.2cm}$$\liminf_{n\to\infty}\frac{1}{m_{l_{n}+1}}\sum_{i=0}^{m_{l_{n}+1}-1}d(x_{i}(u),x_{i}(v))
\geq\lim_{n\to\infty}\frac{m_{l_{n}+1}-m_{l_{n}}}{m_{l_{n}+1}}2\delta=2\delta,                                              \vspace{-0.2cm}$$
which contradicts to (4.8). Hence, $y_u\neq y_{v}$. Consequently, $D_0$ is an uncountable set since $g(\mathcal{E})$ is uncountable.

{\bf Step $2$}.  $D_0$ is a distributionally $\delta$-scrambled set.

Let $u$, $v$, $\{x_i(u)\}_{i=0}^{\infty}$, $\{x_i(v)\}_{i=0}^{\infty}$,
$y_u$, and $y_v$ be specified as those in Step 1.
By applying Lemma 2.12 to (4.7), there exist $E_1,E_2\subset \mathbf{N}$ of density zero such that
\vspace{-0.2cm}$$\lim_{\scriptstyle i\to\infty \atop \scriptstyle i\notin E_1}d(f_{0}^{i}(y_u),x_{i}(u))=0,\;
\;\lim\limits_{\scriptstyle i\to\infty \atop \scriptstyle i\notin E_2}d(f_{0}^{i}(y_{v}),x_{i}(v))=0.\vspace{-0.2cm}$$
Let $\tilde{E}=E_1\cup E_2$. Then $d(\tilde{E})=0$ and
\vspace{-0.2cm}$$\lim_{\scriptstyle i\to\infty \atop \scriptstyle i\notin \tilde{E}}d(f_{0}^{i}(y_u),x_{i}(u))=0,\;\;
\lim\limits_{\scriptstyle i\to\infty \atop \scriptstyle i\notin \tilde{E}}d(f_{0}^{i}(y_{v}),x_{i}(v))=0.                          \eqno(4.10)\vspace{-0.2cm}$$

Firstly, we shall show that $F_{y_u,y_{v}}(\delta)=0$.
By (4.10) there exists an integer $N>0$ such that for each $i>N$ with $i\notin \tilde{E}$,
\vspace{-0.2cm}$$d(f_{0}^{i}(y_{u}),x_{i}(u))<\frac{\delta}{2},\;\;d(f_{0}^{i}(y_{v}),x_{i}(v))<\frac{\delta}{2}.                   \eqno(4.11)\vspace{-0.2cm}$$
By (4.9) and (4.11) one has that for each $N<m_{l_n}\leq i< m_{l_{n}+1}$ with $i\notin \tilde{E}$,
\vspace{-0.2cm}$$d(f_{0}^{i}(y_u),f_{0}^{i}(y_{v}))\geq d(x_i(u),x_i(v))-d(f_{0}^{i}(y_u),x_i(u))-d(f_{0}^{i}(y_{v}),x_i(v))>2\delta-\delta=\delta.\vspace{-0.2cm}$$
Hence, we obtain that
\vspace{-0.2cm}$$\begin{array}{llll}&F_{y_u,y_{v}}(\delta)=\liminf\limits_{n\to\infty}\frac{1}{n}
\sum_{i=0}^{n-1}\chi_{[0,\delta)}(d(f_{0}^{i}(y_u), f_{0}^{i}(y_{v})))\\[1.5ex]
\leq& \liminf\limits_{n\to\infty}\frac{1}{m_{l_{n}+1}}
\sum_{i=0}^{m_{l_{n}+1}-1}\chi_{[0,\delta)}(d(f_{0}^{i}(y_u), f_{0}^{i}(y_{v})))\\[1.5ex]
\leq &\lim\limits_{n\to\infty}\frac{m_{l_{n}}}{m_{l_{n}+1}}+d(\tilde{E})
=0.
\end{array}                                                                     \eqno(4.12)\vspace{-0.2cm}$$

Secondly, we shall show that $F^{*}_{y_u,y_{v}}(\epsilon)=1$ for any $\epsilon>0$. By (4.3) and (4.4)
there exists an increasing sequence $\{s_n\}_{n=0}^{\infty}\subset\mathbf{Z^{+}}$ such that
\vspace{-0.2cm}$$u_{s_n}= v_{s_n},\;n\geq0.                                                                                        \eqno(4.13)\vspace{-0.2cm}$$
It follows from (4.10) that for any $\epsilon>0$, there exists an integer $K>0$ such that for each $i>K$ with $i\notin \tilde{E}$,
\vspace{-0.2cm}$$d(f_{0}^{i}(y_{u}),x_{i}(u))<\frac{\epsilon}{2},\;\;d(f_{0}^{i}(y_{v}),x_{i}(v))<\frac{\epsilon}{2}.               \eqno(4.14)\vspace{-0.2cm}$$
In addition, by (4.5) and (4.13) one has that for each $m_{s_n}\leq i< m_{s_{n}+1}$,
\vspace{-0.2cm}$$d(x_i(u),x_i(v))=d(f_{0}^{i}(u_{s_n}),f_{0}^{i}(v_{s_n}))=0.            \eqno(4.15)\vspace{-0.2cm}$$
It follows from (4.14) and (4.15) that for each $K<m_{s_n}\leq i< m_{s_{n}+1}$ with $i\notin \tilde{E}$,
\vspace{-0.2cm}$$d(f_{0}^{i}(y_u),f_{0}^{i}(y_{v}))\leq d(f_{0}^{i}(y_{u}),x_i(u))+d(x_i(u),x_i(v))+d(x_i(v),f_{0}^{i}(y_{v}))<\epsilon,\vspace{-0.2cm}$$
which implies that
\vspace{-0.2cm}
$$\begin{array}{llll}&F^{*}_{y_u,y_{v}}(\epsilon)=\limsup\limits_{n\to\infty}\frac{1}{n}
\sum_{i=0}^{n-1}\chi_{[0,\epsilon)}(d(f_{0}^{i}(y_u), f_{0}^{i}(y_{v})))\\[1.5ex]
\geq&\limsup\limits_{i\to\infty}\frac{1}{m_{s_{n}+1}}\sum_{i=0}^{m_{s_{n}+1}-1}\chi_{[0,\epsilon)}(d(f_{0}^{i}(y_u), f_{0}^{i}(y_{v})))\\[1.5ex]
\geq&1-\lim\limits_{n\to\infty}\frac{m_{s_{n}}}{m_{s_{n}+1}}-d(\tilde{E})
=1.
\end{array}                                                                    \eqno(4.16)\vspace{-0.2cm}$$
Hence, $D_0$ is a distributionally $\delta$-scrambled set of system (1.1) by (4.12) and (4.16).

Therefore, $D_0$ is an uncountable distributionally $\delta$-scrambled set of system (1.1).
Consequently, system (1.1) is distributionally $\delta$-chaotic.
The proof is complete.\medskip

The following result is a criterion of distributional $\delta$-chaos
induced by some expanding condition for system (1.1).\medskip

\noindent{\bf Theorem 4.6.} {\it Let $(X,d)$ be a compact metric space and $f_n: X\to X$ be continuous for each $n\geq0$.
Assume that there exist two sequences $\{A_n\}_{n=1}^{\infty}$ and $\{B_n\}_{n=1}^{\infty}$ of nonempty closed sets of $X$
satisfying all the assumptions in Lemma {\rm 4.1}, respectively, $A_k\cap B_k=\emptyset$ for some $k\geq1$, and
\vspace{-0.2cm}$$A_{n+1}\cup B_{n+1}\subset f_{n-1}(A_n)\cap f_{n-1}(B_n),\;n\geq1.                                                 \eqno(4.17)\vspace{-0.2cm}$$
Then system {\rm(1.1)} is distributionally $\delta$-chaotic for some $\delta>0$.}\medskip

\noindent{\bf Proof.} It follows from Lemma 4.1 that both $\bigcap_{n=1}^{\infty}A_n$ and $\bigcap_{n=1}^{\infty}B_n$ are singleton sets.
Denote $\bigcap_{n=1}^{\infty}A_n=:\{a\}$ and $\bigcap_{n=1}^{\infty}B_n=:\{b\}$. Then $a\neq b$ since $A_k\cap B_k=\emptyset$ for some $k\geq1$. Let $\delta:=d(a,b)/2$. By the assumption that $d(A_n)\to 0$ and $d(B_n)\to 0$ as $n\to\infty$, there exists an integer $K>0$ such that
\vspace{-0.2cm}$$d(A_{k},B_{k})\geq\delta, \; k>K.                                                                         \eqno(4.18)\vspace{-0.2cm}$$

Let $E\subset\Sigma_{2}^{+}$ be the set satisfying the property in Lemma 2.13. Set $m_1=1$, and
\vspace{-0.2cm}$$m_{n+1}=(2^{n}+1)m_n,\;n\geq1.\vspace{-0.2cm}$$
Denote
\vspace{-0.2cm}$$\Gamma=\big\{\{C_i\}_{i=1}^{\infty}: C_i\in \{A_i, B_i\}, i\geq1\big\}.                                                       \vspace{-0.2cm}$$
Define a map $\varphi: E\to\Gamma$ by $\varphi(\alpha)=\{C_{i}\}_{i=1}^{\infty}$, where $\alpha=(a_0,a_1,\cdots)\in E$,
$$\vspace{-0.2cm}
C_{1}=\begin{cases}
A_1 & \text{ if } a_0=0, \\
B_1 & \text{ if } a_0=1,
\end{cases}$$\vspace{-0.1cm}
and for each $n\geq1$, when $m_{n}<j\leq m_{n+1}$,
$$\vspace{-0.2cm}
C_{j}=\begin{cases}
A_j & \text{ if } a_{n}=0, \\
B_j & \text{ if } a_{n}=1.
\end{cases}\eqno(4.19)$$\vspace{-0.1cm}
It is evident that $\varphi$ is injective. For any $\mathcal{C}=\{C_i\}_{i=1}^{\infty}\in \varphi(E)$, set
\vspace{-0.2cm}$$D_n=\bigcap_{k=0}^{n}f_{0}^{-k}(C_{k+1}).                                                                          \eqno(4.20)\vspace{-0.2cm}$$
It follows from (4.17) that
\vspace{-0.2cm}$$C_{k+1}\subset f_{k-1}(C_k),\;k\geq1.                                                                              \eqno(4.21)\vspace{-0.2cm}$$
By (4.20), (4.21), and the continuity of $f_n$ in compact metric space $(X,d)$ for each $n\geq0$, $D_n$ is a nonempty bounded
and closed set, and satisfies that $D_{n+1}\subset D_{n}$ for each $n\geq1$. Thus, $\bigcap_{n=1}^{\infty}D_n\neq\emptyset$.
Fix one point $x_{\mathcal{C}}\in\bigcap_{n=1}^{\infty}D_n$. This, together with (4.20), yields that
\vspace{-0.2cm}
$$f_{0}^{k}(x_{\mathcal{C}})\in C_{k+1},\; k\geq0.                             \eqno(4.22)\vspace{-0.2cm}$$

Further, we claim that $\mathcal{C}\neq \mathcal{C}'$ if and only if $x_{\mathcal{C}}\neq x_{\mathcal{C}'}$.
It suffices to show the necessity.
Fix any $\mathcal{C}=\{C_i\}_{i=1}^{\infty}, \mathcal{C}'=\{C'_i\}_{i=1}^{\infty}\in \varphi(E)$
with $\mathcal{C}\neq \mathcal{C}'$. Then, there exist $\alpha=(a_0,a_1,\cdots),\alpha'=(a'_0,a'_1,\cdots)\in E$ with $\alpha\neq\alpha'$
such that $\varphi(\alpha)=\mathcal{C}$ and $\varphi(\alpha')=\mathcal{C}'$. So by Lemma 2.13 there exists $n_0>0$ such that
$a_{n_0}\neq a'_{n_0}$ and $m_{n_0}>K$. This, together with (4.18) and (4.19), implies that
$C_{j+1}\neq C'_{j+1}$, and thus $d(C_{j+1},C'_{j+1})\geq \delta$ for each $j$ with $m_{n_0}<j+1\leq m_{n_{0}+1}$.
In addition, it follows from (4.22) that $f_{0}^{j}(x_{\mathcal{C}})\in C_{j+1}$ and $f_{0}^{j}(x_{\mathcal{C}'})\in C'_{j+1}$
for each $j\geq0$. Thus, $x_{\mathcal{C}}\neq x_{\mathcal{C}'}$. Hence, the necessity holds, and then the claim holds.

Denote
\vspace{-0.2cm}$$\tilde{D}=\{x_{\mathcal{C}}: \mathcal{C}\in\varphi(E)\}.                                                                \vspace{-0.2cm}$$
Due to the fact that $E$ is uncountable and the map $\varphi$ is injective, $\varphi(E)$ is uncountable, and
consequently, $\tilde{D}$ is uncountable.

Next, it is to show that $\tilde{D}$ is a distributionally $\delta$-scrambled set.
For any $x,y\in \tilde{D}$ with $x\neq y$, there exist two different points $\mathcal{C}=\{C_i\}_{i=1}^{\infty}, \mathcal{S}=\{S_i\}_{i=1}^{\infty}\in\varphi(E)$ such that $x=x_{\mathcal{C}}$ and $y=x_{\mathcal{S}}$.
Then there exist two different points $\alpha=(a_0,a_1,\cdots), \beta=(b_0,b_1,\cdots)\in E$
such that $\varphi(\alpha)=\mathcal{C}$ and $\varphi(\beta)=\mathcal{S}$. Since $\alpha\neq\beta$, by Lemma 2.13 there exist two sequences
$\{n_i\}_{i=1}^{\infty}, \{l_i\}_{i=1}^{\infty}\subset\mathbf{Z^{+}}$ such that $a_{n_i}\neq b_{n_i}$ and $a_{l_i}=b_{l_i}$ for each $i\geq1$.

Firstly, we shall show that $F_{x,y}(\delta)=0$. Using (4.19) and the fact that $a_{n_i}\neq b_{n_i}$ for each $i\geq1$, one gets that
\vspace{-0.2cm}$$C_{j+1}\neq S_{j+1},\;\;m_{n_i}<j+1\leq m_{n_{i}+1},\vspace{-0.2cm}$$
which, together with (4.18), implies that
\vspace{-0.2cm}
$$d(C_{j+1}, S_{j+1})=d(A_{j+1},B_{j+1})\geq\delta, \;K<m_{n_i}<j+1\leq m_{n_{i}+1}.    \eqno(4.23)\vspace{-0.2cm}$$
It follows from (4.22) and (4.23) that for each $K<m_{n_i}<j+1\leq m_{n_{i}+1}$,
\vspace{-0.2cm}$$d(f_{0}^{j}(x),f_{0}^{j}(y))\geq \delta,                                                                              \vspace{-0.2cm}$$
which results in
\vspace{-0.2cm}$$\liminf_{i\to\infty}\frac{1}{m_{n_{i}+1}}\sum_{j=1}^{m_{n_{i}+1}}\chi_{[0,\delta)}
\bigl(d(f_{0}^{j}(x),f_{0}^{j}(y))\bigr)
\leq\liminf_{i\to\infty}\frac{m_{n_{i}}}{m_{n_{i}+1}}=0.                                                                       \vspace{-0.2cm}$$
Hence, we get that
\vspace{-0.2cm}$$F_{x,y}(\delta)=0.                                                                                                 \eqno(4.24)\vspace{-0.2cm}$$

Secondly, we shall show that $F_{xy}^{*}(\epsilon)=1$ for any $\epsilon>0$. Using (4.19) and the fact that
$a_{l_i}=b_{l_i}$ for each $i\geq1$, one gets that
\vspace{-0.2cm}$$C_{j+1}=S_{j+1},\;\;m_{l_{i}}<j+1\leq m_{l_{i}+1},\;\;i\ge 1.           \eqno(4.25)\vspace{-0.2cm}$$
Because $d(A_j)\to0$ and $d(B_j)\to0$ as $j\to \infty$,
for any $\epsilon>0$, there exists an integer $N>0$ such that
\vspace{-0.2cm}$$d(A_{j})<\epsilon,\; \;d(B_{j})<\epsilon,\;\; j\geq N.                                                      \eqno(4.26)\vspace{-0.2cm}$$
It follows from (4.22), (4.25), and (4.26) that for each $j$ with $N<m_{l_{i}}<j+1\leq m_{l_{i}+1}$,
\vspace{-0.2cm}$$d(f_{0}^{j}(x),f_{0}^{j}(y))<\epsilon,                                                                                 \vspace{-0.2cm}$$
which yields that
\vspace{-0.2cm}$$\limsup_{i\to\infty}\frac{1}{m_{l_{i}+1}}\sum_{j=1}^{m_{l_{i}+1}}\chi_{[0,\epsilon)}
\bigl(d(f_{0}^{j}(x),f_{0}^{j}(y))\bigr)\geq
\limsup_{i\to\infty}\frac{m_{l_{i}+1}-m_{l_{i}}}{m_{l_{i}+1}}=1.                                                              \vspace{-0.2cm}$$
Hence, we obtain that
\vspace{-0.2cm}$$F_{x,y}^{*}(\epsilon)=1.                                                                                           \eqno(4.27)\vspace{-0.2cm}$$
Thus, $\tilde{D}$ is a distributionally $\delta$-scrambled set by (4.24) and (4.27).
Therefore, system (1.1) is distributionally $\delta$-chaotic.
This completes the proof. \medskip

\noindent{\bf Remark 4.7.} The methods used in the proofs of Theorems 4.4 and 4.6 are motivated by
those used in the proofs of Theorem 3.2 in [37] and
Theorem 1 in [35], respectively.\medskip

\bigskip

\noindent{\bf 5. A criterion of distributional chaos in general metric spaces}\bigskip

In this section, we shall establish a criterion of distributional chaos in a sequence for system (1.1)
in a general metric space,
which is induced by a Xiong-chaotic set .\medskip

The following result generalizes Lemma 2 in [42] for ADSs to NDSs.\medskip

\noindent{\bf Lemma 5.1.} {\it Let $P,Q\subset\mathbf{N}$ be two increasing sequences.
Then there exists an increasing subsequence $T$ of $P\cup Q$ such that
\vspace{-0.1cm}$$AR(f_{0,\infty},P)\cap DR(f_{0,\infty},Q)\subset DSR(f_{0,\infty},T).\vspace{-0.5cm}$$}\medskip

\noindent{\bf Proof.} Let $P=\{p_i\}_{i=1}^{\infty}$ and $Q=\{q_i\}_{i=1}^{\infty}$.
The proof is divided into the following two cases.

{\bf Case 1.} Suppose that $P\cap Q=\{r_i\}_{i=1}^{\infty}$ is an infinite sequence. Then one has that
\vspace{-0.2cm}$$AR(f_{0,\infty},P)\cap DR(f_{0,\infty},Q)\subset AR(f_{0,\infty},\{r_i\}_{i=1}^{\infty})\cap DR(f_{0,\infty},\{r_i\}_{i=1}^{\infty})=\emptyset.\vspace{-0.2cm}$$
Clearly, the conclusion holds in this case.

{\bf Case 2.} Suppose that $P\cap Q$ is finite. Without loss of generality, we suppose that
$P\cap Q=\emptyset$. Let
\vspace{-0.2cm}$$n_1=1,\;\; n_{k+1}=2kn_{k},\;\;k\geq1.\vspace{-0.2cm}$$
Now, we define an increasing sequence $T=\{t_i\}_{i=1}^{\infty}$ by $t_1=p_1$ and
\vspace{-0.2cm}$$\{t_i: n_{2k-1}<i\leq n_{2k}\}\subset\{p_j\}_{j=1}^{\infty},\;\;
\{t_i: n_{2k}<i\leq n_{2k+1}\}\subset\{q_j\}_{j=1}^{\infty},\; k\geq1.\vspace{-0.2cm}$$
Fix any $(x,y)\in AR(f_{0,\infty},P)\cap DR(f_{0,\infty},Q)$.
Then, for any $\epsilon>0$ and some $\delta>0$, there exists an integer $N>0$ such that
\vspace{-0.2cm}$$d(f_{0}^{p_i}(x),f_{0}^{p_i}(y))<\epsilon,\;\;d(f_{0}^{q_i}(x),f_{0}^{q_i}(y))>\delta, \; i>N.\vspace{-0.2cm}$$
Thus, for all sufficiently large $k$, one has that
\vspace{-0.2cm}$$d(f_{0}^{t_i}(x),f_{0}^{t_i}(y))<\epsilon,\;\;n_{2k-1}<i\leq n_{2k},\vspace{-0.2cm}$$
which yields that
\vspace{-0.2cm}$$1\geq\limsup_{k\to\infty}\frac{1}{n_{2k}}\sum_{i=1}^{n_{2k}}\chi_{[0,\epsilon)}
(d(f_{0}^{t_i}(x),f_{0}^{t_i}(y)))\geq\lim_{k\to\infty}\frac{n_{2k}-n_{2k-1}}{n_{2k}}=1.\vspace{-0.2cm}$$
Hence, $F_{x,y}^{*}(\epsilon, T)=1$. In addition, for all sufficiently large $k$, we get that
\vspace{-0.2cm}$$d(f_{0}^{t_i}(x),f_{0}^{t_i}(y))>\delta,\;\;n_{2k}<i\leq n_{2k+1},\vspace{-0.2cm}$$
which results in
\vspace{-0.2cm}$$\liminf_{k\to\infty}\frac{1}{n_{2k+1}}\sum_{i=1}^{n_{2k+1}}\chi_{[0,\delta)}
(d(f_{0}^{t_i}(x),f_{0}^{t_i}(y)))\leq\lim_{k\to\infty}\frac{n_{2k}}{n_{2k+1}}=0.\vspace{-0.2cm}$$
Hence, $F_{x,y}(\delta, T)=0$. It follows that $(x,y)\in DSR(f_{0,\infty},T)$.
Therefore, $AR(f_{0,\infty},P)\cap DR(f_{0,\infty},Q)\subset DSR(f_{0,\infty},T)$.
The proof is complete.\medskip

\noindent{\bf Theorem 5.2.} {\it Let $f_n: X\to X$ be continuous for each $n\geq0$. If system {\rm(1.1)}
has an uncountable Xiong-chaotic set with respect to a sequence $P$, then it
is distributionally chaotic in a subsequence of $P$.}\medskip

\noindent{\bf Proof.} Suppose that $S$ is the uncountable Xiong-chaotic set with respect to $P=\{p_i\}_{i=1}^{\infty}$.
Fix two different points $x,y\in S$. Let $F_1: S\to X$ be a constant mapping. Then
there exists a subsequence $Q=\{q_i\}_{i=1}^{\infty}$ of $P$ such that
\vspace{-0.2cm}$$\lim_{i\to\infty}f_{0}^{q_i}(x)=F_1(x)=F_1(y)=\lim_{i\to\infty}f_{0}^{q_i}(y).\vspace{-0.2cm}$$
Then $\lim_{i\to\infty}d(f_{0}^{q_i}(x),f_{0}^{q_i}(y))=0,$
which yields that
\vspace{-0.2cm}$$(x,y)\in AR(f_{0,\infty},Q).                                   \eqno(5.1)\vspace{-0.2cm}$$
Let $F_2: S\to X$ be the identical mapping. Then there exists a subsequence $M=\{m_i\}_{i=1}^{\infty}$ of
$P$ such that
\vspace{-0.2cm}$$\lim_{i\to\infty}f_{0}^{m_i}(x)=F_2(x)=x,\;\;\lim_{i\to\infty}f_{0}^{m_i}(y)=F_2(y)=y.\vspace{-0.2cm}$$
Hence, $\lim_{i\to\infty}d(f_{0}^{m_i}(x),f_{0}^{m_i}(y))>0,$ which implies that
\vspace{-0.2cm}$$(x,y)\in DR(f_{0,\infty},M).                                      \eqno(5.2)\vspace{-0.2cm}$$
It follows from (5.1) and (5.2) that
\vspace{-0.2cm}$$(x,y)\in AR(f_{0,\infty},Q)\cap DR(f_{0,\infty},M).                 \vspace{-0.2cm}$$
Noting that $Q$ and $M$ are independent of $x$ and $y$, one gets that
\vspace{-0.2cm}$$(S\times S)\setminus\triangle\subset AR(f_{0,\infty},Q)\cap DR(f_{0,\infty},M).\vspace{-0.2cm}$$
By Lemma 5.1 there exists an increasing subsequence $T$ of $Q\cup M$ such that
\vspace{-0.2cm}$$(S\times S)\setminus\triangle\subset AR(f_{0,\infty},Q)\cap DR(f_{0,\infty},M)\subset DSR(f_{0,\infty},T).\vspace{-0.2cm}$$
Therefore, $S$ is an uncountable distributionally scrambled set in the sequence $T$,
and consequently, system {\rm(1.1)} is distributionally chaotic in $T$.
This completes the proof.\medskip

\bigskip

\noindent{\bf Acknowledgments}\medskip

This research was supported by the NNSF of China (Grant 11571202).
The first author in the present paper is also supported by China Scholarship Council (File No. 201606220129).

\end{document}